\documentclass[11pt,a4paper]{article}
\usepackage[utf8]{inputenc}
\usepackage[T1]{fontenc}
\usepackage{xcolor}
\usepackage{graphicx}
\usepackage{enumerate}

\usepackage{amsthm,hyperref,lmodern,mathrsfs,amssymb,amsmath,amsfonts}

\usepackage{comment}

\usepackage{geometry}
\newtheorem{thm}{Theorem}

\newtheorem{lemm}{Lemma}

\newtheorem{prop}{Proposition}
\newtheorem{corol}{Corollary}
\newenvironment{dem}{\textbf{Proof}}{\qed }

\newcommand{\R}{\mathbb{R}}
\newcommand{\Per}{\mathrm{Per}}
\newcommand{\E}{\mathbb{E}}

\newcommand{\C}{\mathcal{C}}
\newcommand{\cL}{\mathcal{L}}
\newcommand{\Div}{\mathrm{div}\,}

\title{Generalized perimeters and gradient estimates}

\author{Jordan Serres \footnote{Institut de Mathématiques de Toulouse, France, jserres@insa-toulouse.fr}}

\begin{document}

\maketitle
\begin{center}
\emph{Dedicated to Patrick Cattiaux, my former teacher,\\
from whom I learned the fundamentals of functional analysis.}
\end{center}
\begin{center}
\textbf{Abstract}
\end{center}
\begin{quotation}
\noindent
We use a variational formulation to define a generalized notion of perimeter, from which we derive abstract isoperimetric Cheeger's inequalities via gradient estimates on solutions of Poisson equations. Our abstract framework unifies many existing results and in particular allows us to recover the $W_1-W^{1,1}$ transport inequality, which strengthens the usual transport-information inequality. Conversely, we also prove that Cheeger's inequality implies certain first order Calderón-Zygmund-type gradient estimates.
\end{quotation}

\section{Introduction}

The aim of this article is to revisit the link between isoperimetric inequalities and gradient estimates on solutions of the Poisson equation. Isoperimetric inequalities consist in bounding by below the perimeter of a set by a function of its volume. This problem has attracted considerable interest since antiquity when the usual notions of perimeter and volume, \textit{i.e}. induced by the Lebesgue measure, are considered. Among the incredibly vast literature on the subject, let us mention Ros' survey \cite{ros2001isoperimetric}. However, it is possible to define the notion of perimeter for any non-negative measure, and there are even several different definitions, all of which boil down to the same thing for regular sets. In this work, we will start from the weak formulation of perimeter as the total variation of indicator functions. We refer the reader to Ziemer's book \cite[Chapter 5]{ziemer} for a rigorous overview on functions of bounded variation. This formulation goes back to the work of Caccioppoli \cite{caccioppoli1, caccioppoli2} and De Giorgi \cite{degiorgi53, degiorgi54}, where the perimeter of a set $\Omega\subset\R^d$ is defined by
\[
\mathrm{Per}(\partial\Omega) = \sup \left\{ \int_\Omega \mathrm{div}\,\phi\,|\,\phi\in\mathcal{C}^1_c\left(\R^d,\R^d\right),\, ||\phi||_\infty\leq 1 \right\}
\]
where $\mathcal{C}^1_c\left(\R^d,\R^d\right)$ denotes the set of compactly supported and continuously differentiable functions on $\R^d$, and the integration is with respect to the Lebesgue measure.
Gradient estimates that we will consider consist in showing that for any function $f$ that is sufficiently regular,
\[
|| \nabla f || \leq C\,||\Delta f||
\]
where $C$ does not depend on $f$, and $||\cdot||$ denotes a norm on some domain $\Omega$ with boundary conditions which depends on the context. The most common gradient estimates come from Schauder estimates in the theory of elliptic PDE, which are estimates involving the $\C^2$-norm of $f$, see e.g. \cite{gilbargtrudinger}. In this work, we will consider this approach for the perimeters of general probability measures, and the Laplace operator will therefore be replaced by another second-order elliptic operator. In that case, such estimates are implied by gradient commutation relations at the level of the semigroup generated by the operator, which are equivalent to curvature conditions, see \cite{BGL} for a complete monograph. Such curvature conditions imply strong isoperimetric inequalities, see \cite{bakry1996levy,10.1214/aop/1022874820}, and also Cattiaux and co-authors \cite{10.1214/EJP.v15-754}, and also allow to invert the usual hierarchy between concentration and isoperimetry \cite{milmanequivpoincarepq}.

In this work, we show how the weaker condition of admitting gradient estimates implies isoperimetric Cheeger-type inequalities. This problem is tackled in a general way in Section \ref{sec:main} by introducing an abstract and generalized notion of perimeter, which leads to Theorem \ref{thm:main}. In the rest of the paper, we analyze the consequences of this abstract theorem in different cases. In particular, we show in Section \ref{sec:transport}, that it implies a transport inequality which bounds the Kantorovich-Wasserstein distance $W_1$ between two probability distributions $\mu$ and $\nu$, by the $W^{1,1}$-norm of the density of $\nu$ with respect to $\mu$, see Theorem \ref{thm:transportgradientineq}. This recovers a result in \cite{guillin2009transportationII}. Next, we show that this $W_1-W^{1,1}$ inequality holds for uniformly log-concave distributions, for semi-log-concave probability distributions under an additional assumption, and also for moment measures with uniformly semi-concave moment map. In Section \ref{sec:lpperimeters}, we consider the case of $L^p$-perimeters, which give a natural weakening of the usual isoperimetric problem. We apply our main Theorem \ref{thm:main} to show that under a Poincaré inequality, all $L^p$-perimeters satisfy a Cheeger inequality for $p\geq 2$, see Section \ref{sec:poincareimpliqueLpCheeger}. Finally, in Section \ref{sec:reverse}, we give a partial reverse by showing that Cheeger's inequality implies $L^q-L^p$-Calderón-Zygmund-type gradient estimates for all $p\geq 2$, with $q$ the conjugate exponent of $p$, see Theorem \ref{thm:reversehierarchy}.

\section{The abstract setting}\label{sec:main}

Let $\mu$ be a probability distribution on $\R^d$, let $(E,||\cdot||_E)$ be a space of functions  \hbox{$g:\R^d\to\R$} equipped with a seminorm $||\cdot||_E$, and let $F$ be a space of weakly differentiable functions \hbox{$f:\R^d\to\R$}, equipped with a seminorm $||\cdot||_F$ which takes the following form
$$||f||_F = ||\nabla f||_{\nabla F} $$
where $ ||\cdot||_{\nabla F}$ is a norm on the space of vector-valued functions $\R^d\to\R^d$. We assume that the set $\C^1_c(\R^d,\R)$ of $\C^1$-compactly-supported functions is dense in $F$. We consider a linear operator $\cL : F \to E$, which vanishes on constant functions, and we assume that there exists a constant $C_{E,F}$ such that 
\begin{equation}\label{def:Coperator}
 \forall g\in E,\,\exists f\in F:\, \cL f=g-\int_{\R^d} g\,d\mu \quad \mathrm{and}\quad ||f||_F \leq C_{E,F}\, \left|\left|g -\int_{\R^d} g\,d\mu\right|\right|_E 
\end{equation} 
Intuitively, the constant $C_{E,F}$ is the operator norm of the inverse operator $$\cL^{-1}: E\cap\left\{g\,:\,\int g\,d\mu=0\right\} \to F.$$
Note that  if the equation $\cL f=g-\int_{\R^d} g\,d\mu$ does not admits a solution, or if $\cL^{-1}$ is unbounded, then $C_{E,F}=\infty$. 
We assume that the linear operator $\cL$ admits a divergence form $S$, \textit{i.e}. we assume that there exists some operator $S$ acting on vector valued functions, taking values in $E$, and satisfying that $$\cL(\phi)=S\left(\nabla\phi\right)$$ for all $\phi\in\C^1_c(\R^d,\R)$.
Given a measurable set $\Omega\subset\R^d$, we define its generalized perimeter as follows
\begin{equation}\label{def:generalizedPerimeter}
\Per_{\mu,F}(\partial\Omega) := \sup \left\{\int_\Omega S(\phi)\,d\mu \,|\,\phi\in\mathcal{C}^1_c\left(\R^d,\R^d\right),\, ||\phi||_{\nabla F} \leq 1 \right\}
\end{equation}
We can define its generalized isoperimetric profile function $I_{\mu,F}$ by
\[
I_{\mu,F}(t) := \underset{\mu(\Omega)=t}{\inf} \,\Per_{\mu,F}(\partial\Omega),
\]
which, of course, includes the isoperimetric profile function for the usual perimeter.
We also define, for any two probability distributions $\nu_1,\nu_2$ on $\R^d$, the quantity 
\begin{equation}\label{def:dE}
\delta_E\left(\nu_1 \left| \nu_2\right)\right. := \sup \left\{\int_{\R^d} g\,d\nu_1 -\int_{\R^d} g\,d\nu_2 \,:\,g\in E,\,\left|\left|g -\int_{\R^d} g\,d\nu_2 \right|\right|_E\leq 1 \right\}
\end{equation}
which quantifies the proximity of $\nu_1$ to $\nu_2$. Note that $\delta_E$ is related to usual integro-differential distances between probability distributions such as the total variation distance and the Wasserstein distance, see Sections \ref{sec:usualperimeter}.
For $\Omega\subset\R^d$ with non-zero mass $\mu(\Omega)>0$, we denote by $\mu_\Omega$ the measure $\mu$ restricted to $\Omega$, \textit{i.e}. for all measurable subsets $A\subset\R^d$, it is given by $\mu_\Omega(A)=\mu(\Omega\cap A)/\mu(\Omega)$.
Our main result is the following generalized Cheeger isoperimetric inequality.
\begin{thm}\label{thm:main}
For all probability distributions $\mu$ on $\R^d$, and for all $\Omega\subset\R^d$ with positive mass $\mu(\Omega)>0$, it holds that
\[
\Per_{\mu,F}(\partial\Omega) \geq C_{E,F}^{-1}\,\,\mu(\Omega)\,\, \delta_E\left(\mu_\Omega \left| \mu \right.\right) .
\]
\end{thm}
\noindent
As an immediate corollary, we have the following lower bound on the generalized isoperimetric profile function $I_{\mu,F}$.
\begin{corol}\label{corol:maincorol}
For all probability distributions $\mu$ on $\R^d$ and for all $t\in (0,1)$,
\[
I_{\mu,F}(t) \geq t\,\, \underset{E}{\sup} \left( C_{E,F}^{-1}\,\underset{\mu(\Omega)=t}{\inf} \delta_E\left(\mu_\Omega \left| \mu \right.\right)\right),
\]
where the supremum runs over all seminormed spaces $(E,||\cdot||_E)$ of functions on $\R^d$, and the constant $C_{E,F}$ denotes the operator norm of $\cL^{-1}$ given in \eqref{def:Coperator}.
\end{corol}
\noindent
Note that if $\cL^{-1}$ is unbounded or undefined on $E$, then the inverse of its operator norm is zero, and therefore the corresponding term in the supremum is not taken in account.\newline\newline
{\textbf{Proof of Theorem \ref{thm:main}}.
We can write
\begin{align*}
\Per_{\mu,F}(\partial\Omega) &= \sup \left\{\int_\Omega S(\phi)\,d\mu \,|\,\phi\in\mathcal{C}_c^1\left(\R^d,\R^d\right),\, ||\phi||_{\nabla F} \leq 1 \right\}\\
&= C_{E,F}^{-1}\, \mu(\Omega)\,\sup \left\{\int_\Omega S(\phi)\,d\mu_\Omega\,|\,\phi\in\mathcal{C}_c^1\left(\R^d,\R^d\right),\, ||\phi||_{\nabla F}\leq C_{E,F} \right\}\\
&\geq C_{E,F}^{-1}\, \mu(\Omega)\,\sup \left\{\int_\Omega S(\nabla \psi) \,d\mu_\Omega\,|\,\psi\in\mathcal{C}_c^1\left(\R^d,\R \right),\, ||\nabla \psi||_{\nabla F}\leq C_{E,F} \right\}\\
&= C_{E,F}^{-1}\, \mu(\Omega)\,\sup \left\{\int_\Omega \cL f \,d\mu_\Omega\,|\,f\in F,\, ||f||_F\leq C_{E,F} \right\}\\
&\geq C_{E,F}^{-1}\, \mu(\Omega)\,\sup \left\{\int_\Omega \cL f\,d\mu_\Omega\,|\,f\in F\,\mathrm{given\, in\, \eqref{def:Coperator}\, for\, some\,} g\in E,\, \left|\left|g - \int_{\R^d} g\,d\mu\right|\right|_E \leq 1 \right\}\\
&= C_{E,F}^{-1}\, \mu(\Omega)\,\sup \left\{\int_\Omega \left(g -\int_{\R^d} g\,d\mu \right) d\mu_\Omega\,|\, g\in E,\, \left|\left|g - \int_{\R^d} g\,d\mu\right|\right|_E \leq 1 \right\}\\
&= C_{E,F}^{-1}\, \mu(\Omega)\,\delta_E\left(\mu_\Omega \left| \mu \right.\right)
\end{align*}
where, on line $4$, we have used the assumption that the set of $\C^1$-compactly supported functions is dense in $(F,||\cdot||_F)$, and on the last line, we have used that $\mu_\Omega$ is a probability distribution.
The proof is complete. $\qed$}

\section{The usual perimeter}\label{sec:usualperimeter}

In this section, we consider the case of a probability distribution $\mu=e^{-V}dx$, $V\in\C^2(\R^d)$, with its usual perimeter $\mu^+$, \textit{i.e}. we consider Definition \eqref{def:generalizedPerimeter} with $(F,||\cdot||_F)$ the set of Lipschitz functions equipped with the Lipschitz seminorm $||f||_F = ||\nabla f||_{\infty}$, and $S(\phi)=\Div\phi - \nabla V\cdot \phi$. Indeed, it is well known that in that case, the quantity 
\[
\mu^+\left(\partial\Omega\right) = \sup \left\{ \int_\Omega \Div \phi-\nabla V\cdot \phi \,d\mu \,:\,\phi\in\mathcal{C}^1_c\left(\R^d,\R^d\right),\, ||\phi||_\infty\leq 1 \right\}
\]
is the total variation of the indicator function $\textbf{1}_\Omega$ with respect to the measure $\mu$, which is actually the $\mu$-perimeter of $\partial\Omega$, see e.g. \cite{maggi2012sets}. In the next two sections, we will consider the cases where $E$ is the set of bounded functions, giving the usual Cheeger inequalities, see Section \ref{sec:cheegerusuel}, and where $E$ is the set of Lipschitz functions, giving the $W_1-W^{1,1}$ transport inequality, see Section \ref{sec:transport}.

\subsection{Cheeger's inequality}\label{sec:cheegerusuel}

When $E=L^\infty$ is the set of bounded functions equipped with the uniform norm $||\cdot||_\infty$, the quantity $\delta_E$ defined by \eqref{def:dE} controls the total variation distance: For all $\Omega$ such that $\mu(\Omega)>0$,
\[
\delta_{L^\infty} \left(\mu_\Omega \left| \mu \right.\right) \geq \frac{1}{2} d_{TV} ( \mu_\Omega,\mu ),
\]
where $$ d_{TV} ( \mu_\Omega,\mu ) = \sup \left\{\int_{\R^d} g\,d\mu_\Omega -\int_{\R^d} g\,d\mu \,:\,\,\left| \left|g \right| \right|_\infty \leq 1 \right\} $$ denotes the total variation distance between the two probability distributions $\mu$ and $\mu_\Omega(\cdot)= (1/\mu(\Omega))\mu(\Omega\cap\cdot)$.
This is an immediate consequence of the fact that $\left|\left|g-\int g\,d\mu\right|\right|_\infty \leq 2\left|\left|g\right|\right|_\infty $. Moreover, by taking the combination of indicators $g=\textbf{1}_\Omega-\textbf{1}_{\R^d\setminus\Omega}$, it is easy to see that $d_{TV}(\mu_\Omega, \mu ) = 2\left(1-\mu(\Omega)\right).$ Therefore, in that case, Corollary \ref{corol:maincorol} gives the usual Cheeger inequality
\[
I_\mu(t) \geq  C_{L^\infty,Lip}^{-1}\,t(1-t)
\]
The fact that gradient estimates imply Cheeger's inequality was already noticed by Wu \cite[Corollary 2.2]{wu2009gradient}, who showed that such gradient estimates hold on a complete connected Riemannian manifold when $V\in\C^2$ admits a second moment, and the constant $C_{L^\infty,Lip}$ depends on the Ricci curvature.
Let us end this section by mentioning that this method is related to Stein's method when $\mu=\gamma$ is the standard one-dimensional normal distribution.
Indeed, in that case, the perimeter is given by
\[\gamma^+(\Omega) =\sup \left\{\int_\Omega v'(x)-x\,v(x)\,d\gamma(x)\,:\,v\in\mathcal{C}^1_c\left(\R,\R\right),\, ||v||_\infty\leq 1 \right\},
\]
and Stein's lemma \cite[Lemma 2.5]{fundamStein} gives that for any bounded function $g$, there is a solution of the Stein equation $v'(x)-xv(x)= g(x) - \int g\,d\gamma$ satisfying $||v||_\infty\leq \sqrt{\frac{\pi}{2}} ||g(x) - \int g\,d\gamma||_\infty $. Therefore $C_{L^\infty,Lip}^{-1}=\sqrt{\frac{2}{\pi}}$ which, up to a factor $1/2$, gives the sharp Cheeger constant for the Gaussian.

\subsection{The $W_1 - W^{1,1}$ transport inequality}\label{sec:transport}

When $E=\mathrm{Lip}$ is the set of Lipschitz functions equiped with the Lipschitz seminorm $||\nabla g||_{\infty}$, the quantity $\delta_E$ is the $L^1$-Wasserstein distance $$\delta_E (\mu_\Omega|\mu) = W_1(\mu_\Omega,\mu), $$ where
\[
W_1(\mu_\Omega,\mu) := \sup \left\{\int_{\R^d} g\,d\mu -\frac{1}{\mu(\Omega)}\int_{\Omega} g\,d\mu \,:\,||\nabla g||_\infty\leq 1 \right\}
\]
and Theorem \ref{thm:main} gives that for all $\Omega\subset\R^d$ with non-zero mass, 
\begin{equation}\label{eq:applithm1pourtransport}
\mu^+(\partial\Omega) \geq C_{Lip,Lip}^{-1}\,\,\mu(\Omega)\,W_1(\mu ,\mu_\Omega)
\end{equation}
from which we deduce the following.
\begin{thm}\label{thm:transportgradientineq}
For all probability distributions $\mu$ on $\R^d$, and for all functions $f:\R^d\to\R_+$ such that $\int f\,d\mu=1$, it holds that
\begin{equation}\label{eq:transportgradientineq}
W_1\left( \mu, f \mu \right) \leq C_{Lip,Lip} \int_{\R^d} \left| \nabla f\right|\,d\mu
\end{equation}
\end{thm}
\noindent
Note that when the constant $C_{Lip,Lip}$ is not bounded, then Inequality \eqref{eq:transportgradientineq} is just empty. 
Let us point out that Theorem \ref{thm:transportgradientineq} recovers \cite[Theorem 3.1]{guillin2009transportationII} of Guillin, Léonard, Wang and Wu.
Note also that Lipschitz-to-Lipschitz gradient estimates are a vast topic in the literature, including Kuwada's duality theorem, stating that for semigroups, they are equivalent to a Wasserstein contraction (see \cite{kuwada2010duality}). We now give the proof of Theorem \ref{thm:transportgradientineq}.\newline\newline
\begin{dem}
Let $h:\R^d\to \R$ be a $1$-Lipschitz function. We can compute that
\begin{align*}
\int_{\R^d} \left| \nabla f\right|\,d\mu & = \int_0^\infty \mu^+\left(f>t \right)\,dt\\
& \geq C_{Lip,Lip}^{-1}\,\int_0^\infty \mu(f>t)\,W_1\left(\mu ,\mu_{\{f>t\}}\right) dt\\
& \geq C_{Lip,Lip}^{-1}\,\int_0^\infty \mu(f>t)\,\int_{\R^d} h(x)\left( 1 - \frac{\textbf{1}_{\{f>t\}}(x)}{\mu(f>t)} \right)d\mu(x)\,dt\\
& = C_{Lip,Lip}^{-1}\,\int_{\R^d} h(x) \int_0^\infty \left( \mu(f>t) - \textbf{1}_{\{f>t\}}(x) \right)dt\,d\mu(x)\\
& = C_{Lip,Lip}^{-1}\,\int_{\R^d} h(x) \left( \int f\,d\mu - f(x) \right) d\mu(x)\\
& = C_{Lip,Lip}^{-1}\,\int_{\R^d} h(x) \left( 1 - f(x) \right) d\mu(x)\\
\end{align*}
where we used the co-area formula on lines 1 and 5, Inequality \eqref{eq:applithm1pourtransport} on line 2, and the fact that $\int f\,d\mu=1$ on the last line. Since this holds for all $1$-Lipschitz functions $h$, by taking the supremum, we deduce that $$ \int_{\R^d} \left| \nabla f\right|\,d\mu \geq C_{Lip,Lip}^{-1}\, W_1\left( \mu, f \mu \right) $$ which is exactly Inequality \ref{eq:transportgradientineq}.
\end{dem}\newline
Notice that the converse statement, that is the fact that the transport Inequality \eqref{eq:transportgradientineq} implies the isoperimetric-type inequality \eqref{eq:applithm1pourtransport} is also true, as can be seen by taking Lipschitz approximations of the indicator function $\textbf{1}_\Omega$ in \eqref{eq:transportgradientineq} and passing to the limit.
Let us also point out that Inequality \eqref{eq:transportgradientineq} was first introduced in \cite{guillin2009transportationII}. We refer the reader to \cite{gozlanleonardotettransportineq} for a survey on transport inequalities, and to \cite{guillin2009transportation,guillin2009transportationII} for comparison between the transport inequality \eqref{eq:transportgradientineq} and others transports inequalities. In particular, let us underline that the $W_1-W^{1,1}$ inequality \eqref{eq:transportgradientineq} can be seen as a reinforcement of the usual transport-information inequality. Indeed, by using the Cauchy-Schwarz inequality, we can write that
\begin{align*}
W_1\left( \mu, f \mu \right) &\leq C_{Lip,Lip} \int_{\R^d} \left| \nabla f\right|\,d\mu\\
& = C_{Lip,Lip} \int_{\R^d} \frac{\left| \nabla f\right|}{\sqrt{f}} \sqrt{f}\,d\mu\\
& \leq C_{Lip,Lip} \sqrt{\int_{\R^d} \frac{\left| \nabla f\right|^2}{f}\,d\mu} \sqrt{\int_{\R^d} f\,d\mu}\\
& = \sqrt{C_{Lip,Lip}^2 \int_{\R^d} \frac{\left| \nabla f\right|^2}{f}\,d\mu}.
\end{align*}
The next subsection of this section will now examine, from the literature, how the constant $C_{Lip, Lip}$ can be explicitly bounded in the case of uniformly log-concave and semi-log-concave probability distributions, the second subsection will show how to derive the transport $W_1$-$W^{1,1}$ inequality for moment measures, and finally at the third subsection, we will show that the transport $W_1$-$W^{1,1}$ inequality implies a Cheeger-type inequality. 

\subsubsection{Application to semi-log-concave distributions}

Let $d\mu(x)=e^{-V(x)}dx$ be a probability distribution, and for all $x\in\R^d$, let $\kappa(x)\in\R$ denote the smallest eigenvalue of the symmetric matrix $\nabla^2 V(x)$, so we have $$\forall x\in\R^d,\quad \nabla^2 V(x) \geq \kappa(x) I_d.$$ Then it is possible to show, see e.g. \cite[Section 4.9.2]{BGL}, that for all Lipschitz functions $g$, the equation
\[
\Delta f -\nabla V\cdot \nabla f = g -\int_{\R^d} g\,d\mu
\]
admits a solution $f$ satisfying 
\[
||\nabla f||_{\infty} \leq ||\nabla g||_{\infty}\,\sup_{x\in\R^d} \int_0^\infty \E \left[ e^{-\int_0^t \kappa(X_s)\,ds} \right] dt,
\]	 
where $(X_s)_{s\geq 0}$ is the Markov diffusion generated by the operator $\cL=\Delta -\nabla V\cdot \nabla $ and starting at $X_0=x$.
Let us mention that this type of probabilistic representation was used by Cattiaux \textit{et al}., to study the long-term behavior of diffusion processes in \cite{cattiaux2022self}.
In particular, this estimate assures us that the constant $C_{Lip,Lip}$ satisfies
\begin{equation}\label{eq:FeynmanKac}
C_{Lip,Lip}\leq  \sup_{x\in\R^d} \int_0^\infty \E \left[ e^{-\int_0^t \kappa(X_s)\,ds} \right] dt
\end{equation}
and Theorem \ref{thm:transportgradientineq} gives that as soon as this supremum is finite, then for all functions $f:\R^d\to\R_+$ such that $\int f\,d\mu=1$, it holds
\[
W_1\left( \mu, f \mu \right) \leq \left( \sup_{x\in\R^d} \int_0^\infty \E \left[ e^{-\int_0^t \kappa(X_s)\,ds} \right] dt\right) \int_{\R^d} \left| \nabla f\right|\,d\mu.
\]
In particular, if for some $\rho>0$, $d\mu(x)=e^{-V(x)}dx$ is $\rho$-uniformly log-concave, \textit{i.e}. $\nabla V \geq \rho I_d$ in the sense of quadratic forms with $I_d$ the identity matrix, then we can take $\kappa(x)=\rho$, and \eqref{eq:FeynmanKac} gives
\[
C_{Lip,Lip}\leq 1/\kappa,
\]
and Theorem \ref{thm:transportgradientineq} gives that for all functions $f:\R^d\to\R_+$ such that $\int f\,d\mu=1$,
\[
W_1\left( \mu, f \mu \right) \leq \frac{1}{\kappa} \int_{\R^d} \left| \nabla f\right|\,d\mu.
\]

\subsubsection{Application to moment measures}

In this section, we suppose that $\mu$ is a moment measure, \textit{i.e}. that for some $\varphi : \R^d\to \R$, it can be written as the pushforward by $\nabla\varphi$ of the probability distribution $d\nu(x):=e^{-\varphi(x)}dx$, 
\[
d\mu(x) = \nabla\varphi^{\#}\left( e^{-\varphi(x)}dx\right),
\]
and $\nabla\varphi$ is then called the moment map of $\mu$. Such a function exists if $\mu$ is centered, admits a finite first moment and is not supported on a hyperplane, and moreover it is the Brenier map from optimal transport theory, see \cite{cordero2015moment}. We have the following.
\begin{thm}
If $\mu$ is a moment measure and $\nabla\varphi$ denotes its moment map, then it satisfies the following transport $W_1 - W^{1,1}$ inequality for all non-negative functions $f: \R^d\to\R_+$ such that $\int f\,d\mu=1$,
\[
W_1\left(\mu,f\mu\right) \leq 2\, \sup_{x\in\R^d} \lambda_{\max}\left(\nabla^2\varphi(x)\right)\, \int_{\R^d} \left| \nabla f\right|\,d\mu,
\]
where $\lambda_{\max}\left(\nabla^2\varphi(x)\right)$ denotes the largest eigenvalue of the matrix $\nabla^2\varphi(x)$.
\end{thm}
\noindent
Note that since $\nabla\varphi$ is the Brenier map from optimal transport theory, it follows that $\varphi$ is convex, and therefore for all $x\in\R^d$, $\lambda_{\max}\left(\nabla^2\varphi(x)\right)\geq 0$. Note also that the standard $d$-dimensional Gaussian distribution is the only fixed point of the moment measure functional, \textit{i.e.} the only distribution such that $\nabla\varphi = I_d$. As a consequence, it is natural to expect that the quantity $I_d - \nabla^2\varphi$ characterizes the proximity between $\mu$ and the standard normal distribution. This observation was put forward by Fathi in \cite{fathi2019stein} where it was shown that $\nabla^2\varphi\left(\nabla\varphi^{-1}\right)$ is a Stein kernel for $\mu$, from which it follows that the $L^2$ norm of $I_d - \nabla^2\varphi$ controls the $L^2$-Wasserstein distance between $\mu$ and the Gaussian. By way of comparison, our bound involving the supremum of $\nabla^2\varphi$ is much more restrictive.\newline\newline
\begin{dem}
By using the integration by parts formula, it is straigtforward to see that there exists some vector-valued function $W:\R^d\to\R^d$ such that the operator
\[
\Lambda f := \langle \nabla^2\varphi^{-1}, \nabla^2 f \rangle_{HS} - W\cdot \nabla f
\]
satisfies for all smooth functions $f,g:\R^d\to\R$,
\[
\int_{\R^d} f\Lambda g\,d\nu = -\int_{\R^d} \nabla^{2}\varphi^{-1}\left(\nabla f,\nabla g \right)d\nu.
\]
Let us also define the operator $S$ which acts on vector-valued functions $\phi:\R^d\to\R^d$ as $S(\phi) = \langle \nabla^2\varphi^{-1}, \nabla \phi \rangle_{HS} - W\cdot \phi$.
It has been shown by Kolesnikov in \cite{Kolesnikov12} that the metric measure space $(\R^d,\nabla^2\varphi,\,\nu)$ satisfies the curvature-dimension condition $CD(1/2,\infty)$. In particular, this implies that the solution of the Poisson equation $\Lambda f = g - \int_{\R^d} g\,d\nu$ satisfies the gradient estimate $||\nabla^2\varphi^{-1}(\nabla f,\nabla f)||_{\infty} \leq 2 ||\nabla^2\varphi^{-1}(\nabla g,\nabla g)||_{\infty} $, see e.g. \cite{BGL}. Therefore we have that the operator norm $C_{Lip_\varphi,Lip_\varphi}$ satisfies $C_{Lip_\varphi,Lip_\varphi} \leq 2 $, where $Lip_\varphi$ denotes the space of Lipschitz functions with respect to the Hessian metric $\nabla^2\varphi$, and Theorem \ref{thm:transportgradientineq} gives the following transport inequality for $\nu=e^{-\varphi}dx$
\begin{equation}\label{eq:transportfornu}
W_{1,\varphi}\left( \nu, f \nu \right) \leq 2 \int_{\R^d} \sqrt{\nabla^{2}\varphi^{-1}\left(\nabla f,\nabla f \right)} \,d\nu
\end{equation}
which holds for all non-negative functions $f: \R^d\to\R_+$ such that $\int f\,d\nu=1$, and where 
\[
W_{1,\varphi}\left(\alpha,\beta\right) := \sup \left\{\int_{\R^d} h\,d\alpha -\int_{\R^d} h\,d\beta \,:\,\left|\left|\nabla^{2}\varphi^{-1}\left(\nabla h,\nabla h \right)\right|\right|_\infty\leq 1  \right\}
\]
denotes the $L^1$-Kantorovich distance with respect to the distance induced by the Hessian metric $\nabla^2\varphi$. Using the fact that $\nabla\varphi$ is the moment map of $\mu$, we can see that for all non-negative functions $f: \R^d\to\R_+$ such that $\int f\,d\mu=1$,
\begin{align*}
W_1\left(\mu,f\mu\right) &= W_{1}\left( \nabla\varphi^{\#}(\nu), f \nabla\varphi^{\#}(\nu)\right)\\
&= \sup \left\{\int_{\R^d} g(\nabla\varphi)\,d\nu -\int_{\R^d} g(\nabla\varphi)f(\nabla\varphi)\,d\nu \,:\,||\nabla g||_\infty\leq 1 \right\} \\
& \leq \sup \left\{\int_{\R^d} h\,d\nu -\int_{\R^d} h f(\nabla\varphi)\,d\nu \,:\,\left|\left|\sqrt{\nabla^{2}\varphi^{-1}\left(\nabla h,\nabla h \right)}\right|\right|_\infty\leq \sup_{x\in\R^d} \sqrt{\lambda_{\max}\left(\nabla^2\varphi(x)\right)}  \right\}\\
& \leq  \sup_{x\in\R^d} \lambda_{\max}\left(\nabla^2\varphi(x)\right)^{1/2}\, W_{1,\varphi}\left(\nu, (f\circ\nabla\varphi)\,\nu \right)\\
& \leq 2 \sup_{x\in\R^d} \lambda_{\max}\left(\nabla^2\varphi(x)\right)^{1/2} \int_{\R^d} \sqrt{\nabla^{2}\varphi^{-1}\left(\nabla (f\circ\nabla\varphi),\nabla (f\circ\nabla\varphi) \right)} \,d\nu\\
& \leq 2 \sup_{x\in\R^d} \lambda_{\max}\left(\nabla^2\varphi(x)\right)^{1/2} \int_{\R^d} \sqrt{\nabla^{2}\varphi \left((\nabla f)\circ\nabla\varphi,(\nabla f)\circ\nabla\varphi \right)} \,d\nu\\
& \leq 2 \sup_{x\in\R^d} \lambda_{\max}\left(\nabla^2\varphi(x)\right) \int_{\R^d} \left| \nabla f\right|\,d\mu, 
\end{align*}
where we used the definition of $W_{1,\varphi}$ on line 4, and Inequality \eqref{eq:transportfornu} on line 5.
The proof is complete.
\end{dem}

\subsubsection{A Cheeger-type isoperimetric inequality}

In this section, we show that the $W_1 - W^{1,1}$ transport inequality of Theorem \ref{thm:transportgradientineq} implies a different version of Cheeger's inequality.
\begin{prop}
If the probability distribution $\mu$ on $\R^d$ satisfies the $W_1 - W^{1,1}$ transport inequality, \textit{i.e}. for all functions $f:\R^d\to\R_+$ such that $\int f\,d\mu=1$,
\[
W_1\left( \mu, f \mu \right) \leq C \int_{\R^d} \left| \nabla f\right|\,d\mu
\]
then $\mu$ satisfies the following version of Cheeger's inequality, for all measurable subsets $\Omega\subset\R^d$,
\[ 
\mu^+(\partial\Omega) \geq  C^{-1}\, \min(\mu(\Omega),1-\mu(\Omega)) \int_{\R^d}d(x,\partial\Omega)\,d\mu(x).
\]
\end{prop} 
\noindent
\begin{dem}
The proof consists in applying the $W_1 - W^{1,1}$ transport inequality to Lipschitz approximations of the indicator function $\mu(\Omega)^{-1}\textbf{1}_\Omega$, and then observing that the distance $W_1$ satisfies
\[
W_1(\mu ,\mu_\Omega) \geq \min \left(1,\frac{1-\mu(\Omega)}{\mu(\Omega)}\right)\, \int_{\R^d}d(x,\partial\Omega)\,d\mu(x).
\]
Indeed, by choosing
\[
g(x) =  \left\{ 
\begin{array}{rcl} 
d(x,\partial\Omega) \quad & \, x\in\Omega \\ 
-d(x,\partial\Omega) \quad & \, x\in\R^d\setminus\Omega \\ 
\end{array} \right.
\]
which is $1$-Lipschitz, one can see that
$$W_1(\mu ,\mu_\Omega) \geq \frac{1-\mu(\Omega)}{\mu(\Omega)}\int_\Omega d(x,\partial\Omega)\,d\mu + \int_{\R^d\setminus\Omega} d(x,\partial\Omega)\,d\mu.$$
Now if $\mu(\Omega)\in (0,1/2)$, we have that $\frac{1-\mu(\Omega)}{\mu(\Omega)}\geq 1$ and therefore $W_1(\mu_\Omega \left| \mu \right.)\geq \int_{\R^d} d(x,\partial\Omega)\,d\mu$, and similarly, if $\mu(\Omega)\in (1/2,1)$, then $1\geq \frac{1-\mu(\Omega)}{\mu(\Omega)}$, so we get that 
$W_1(\mu_\Omega \left| \mu \right.)\geq \frac{1-\mu(\Omega)}{\mu(\Omega)} \int_{\R^d}d(x,\partial\Omega)\,d\mu$, achieving the proof.
\end{dem}\newline\newline
It is worth mentioning that Corollary \ref{corol:maincorol} gives the following bound on the isoperimetric profile function
\[
I_{\mu}(t) \geq C_{Lip, Lip}^{-1}\,\min \left(t,1-t\right) \underset{\mu(\Omega)=t}{\inf} \int_{\R^d}d(x,\partial\Omega)\,d\mu.
\]
Let us make a few comments on the minimization problem
\begin{equation}\label{eq:medianset}
\underset{\mu(\Omega)=t}{\inf} \int_{\R^d}d(x,\partial\Omega)\,d\mu
\end{equation}
Firstly, if $\partial\Omega$ is not regular enough, the infimum may be zero. For example, let us take the uniform distribution on $[0,1]$, and a dyadic decomposition whose volume remains constant but whose diameter of connected components becomes zero. This Cheeger-type inequality therefore only makes sense if we restrict ourselves to sufficiently regular sets, e.g. at least connected. Secondly, note that if we take $t=0$ and restrict the infimum to singleton sets, then the minimizers of Problem \eqref{eq:medianset} are exactly the medians of $\mu$. Consequently, even in the case $t>0$, the minimizers among the connected sets can be interpreted as $t$-mass medians, and the infimum itself as the $L^1$-variance for $t$-mass sets. In particular, for the multivariate standard normal distribution, the minimizer will be the ball centered at zero whose Gaussian mass is equal to $t$.

\section{The $L^p$-perimeters, $p\geq 1$}\label{sec:lpperimeters}

The usual perimeter is defined as the total variation of indicator functions, \textit{i.e}. as the operator norm of $\phi\mapsto\int_\Omega \Div \phi$ on the space $L^\infty$. In this section, we define the $L^p$-perimeter $\mu^{+,p}$, $p\geq 1$, of a probability distribution $d\mu(x)=e^{-V(x)}dx$, as the $p$-total variation of the gradient of indicator functions
\[
\mu^{+,p} (\partial\Omega) := \sup \left\{\int_\Omega \Div\phi-\nabla V\cdot \phi\,d\mu\,|\,\phi\in\mathcal{C}_c^1\left(\R^d,\R^d\right),\, ||\phi||_{L^q(\mu)}\leq 1 \right\}
\]
where $q$ is the conjugate exponent of $p$, \textit{i.e}. $1/p+1/q=1$. In the terminology of the abstract setting of Section \ref{sec:main}, this means that we consider Definition \eqref{def:generalizedPerimeter} with $F$ equals to the Sobolev space $W^{1,q}(\mu)$ equipped with its seminorm $||f||_F = ||\nabla f||_{L^q(\mu)}$, and the divergence form $S$ is given by 
\[
S(\phi) = \Div \phi - \nabla V\cdot \phi.
\]
Note that from this definition, we have that $\mu^{+,p} (\partial\Omega) = \left(\int_{\R^d} |\nabla \mathbf{1}_\Omega|^p d\mu\right)^{1/p}$, and therefore the usual perimeter is the $L^1$-perimeter $\mu^+=\mu^{+,1}$. Note moreover that since $\mu$ is a probability distribution, we have the inclusions $L^r(\mu)\subset L^s(\mu)$, for $1\leq s \leq r \leq \infty$, and therefore
\[
\forall \Omega,\quad \mu^{+,r}(\partial \Omega) \geq \mu^{+,s} (\partial \Omega).
\]
In other words, lower bounds on the usual perimeter imply lower bounds for all other perimeters, so that isoperimetric inequalities concerning $L^p$-perimeters, $p>1$, represent a weakened version of the usual isoperimetric problem.

\subsection{$L^p$-Cheeger's inequalities}\label{sec:LpCheeger}

In this section, we consider the case of a probability distribution $\mu=e^{-V}dx$, $V\in\C^2(\R^d)$, with its $p$-perimeter $\mu^{+,p}$, $p\geq 1$. By using the notation of the abstract setting of Section \ref{sec:main}, we take $E=L^p(\mu)$ equipped with the $L^p$-norm. The gradient estimate given in \eqref{def:Coperator} then becomes the following first order Calderón-Zygmund estimate,
\begin{equation}\label{eq:CZgeneral}
||\nabla f||_{L^q(\mu)} \leq C_{q,p} ||\cL f||_{L^p(\mu)}
\end{equation}
where $C_q$ is a constant independant of $f$, and $\cL f = S(\nabla f) = \Delta f -\nabla f\cdot \nabla V$. Note that our terminology comes from the usual Calderón-Zygmund estimate which involves the Hessian of $f$ instead of its gradient. Theorem \ref{thm:main} then gives
\[
\mu^{+,p}(\partial\Omega) \geq C_{q,p}^{-1}\,\, \mu(\Omega)\,\, \delta_{L^p(\mu)} \left(\mu_\Omega \left| \mu \right. \right)
\]
where $\delta_{L^p(\mu)}$ is defined in \eqref{def:dE} with $E=L^p(\mu)$.
Since $||\textbf{1}_\Omega - \mu(\Omega)||_{L^p(\mu)}\leq 1$, we can take $g=\textbf{1}_\Omega$ in Definition \ref{def:dE}, and we obtain $\delta_{L^q}(\mu_\Omega \left| \mu \right.) \geq (1-\mu(\Omega))$, from which we deduce the following $p$-Cheeger inequality
\begin{equation}\label{eq:LpCheeger}
\mu^{+,p}(\partial\Omega) \geq  C_{q,p}^{-1} \, \mu(\Omega)\left( 1-\mu(\Omega)\right)
\end{equation}
Note that to obtain Cheeger's inequality for the usual perimeter with this methodology, we would need a first order Calderón-Zygmund-type estimate with $p=1$ and $q=\infty$, which is a stronger gradient estimate than the one used in Section \ref{sec:cheegerusuel}.

\subsection{First order Calderón-Zygmund estimates}\label{sec:CZ}

We say that a probability distribution $\mu$ satisfies the Poincaré inequality with constant $\lambda_2$ if for all $f\in H^1(\mu)$, it holds
\begin{equation}\label{def:Poincare}
\lambda_2\,\int \left( f - \int f\,d\mu \right)d\mu \leq \int \left|\nabla f\right|^2 d\mu.
\end{equation}
It is well known that the Poincaré inequality is equivalent to the following convergence in $L^2$,
\[
\left|\left| P_t(f) - \int f\,d\mu \right|\right|_{L^2(\mu)} \leq e^{-\lambda_2\,t} \left|\left| f - \int f\,d\mu \right|\right|_{L^2(\mu)}
\]
where $P_t$ is the semigroup generated by the operator $\cL$ for which $\mu$ is a reversible probability distribution, see e.g. \cite{ane2000inegalites} for details. Moreover, Cattiaux and his collaborators showed in \cite{CGR10} that the Poincaré inequality also implies the following $L^p$ convergence of semigroups
\begin{equation}\label{eq:cvLpsemigrp}
\left|\left| P_t(f) - \int f\,d\mu \right|\right|_{L^p(\mu)} \leq K_p\, e^{-\lambda_p\,t} \left|\left| f - \int f\,d\mu \right|\right|_{L^p(\mu)}
\end{equation}
where $K_p=2^{1-r_p}$, $\lambda_p=r_p\lambda_2$ and $r_p=2\min\left(\frac{1}{p},1-\frac{1}{p}\right) $. In the following lemma, we show that the Poincaré inequality implies the first order Calderón-Zygmund estimate for all $p\geq 2$.
\begin{lemm}\label{lem:CZforL}
Let $\mu=e^{-V}dx$, $V\in\C^2(\R^d)$, be a probability distribution, let $p\geq 2$, and let $g\in L^q(\mu)$, for $q\in[1,2]$ the conjugate exponent of $p$. If $\mu$ satisfies a Poincaré inequality with constant $\lambda_2$, then the Poisson equation $$\cL f = g - \int g\,d\mu, $$ admits a solution $f$ satisfying the following first order Calderón-Zygmund estimate
\begin{equation}\label{eq:lemCZforL}
||\nabla f||_{L^q(\mu)} \leq C_{p}\, \left|\left|g - \int g\,d\mu \right|\right|_{L^p(\mu)}
\end{equation}
where $C_{p} = \sqrt{2}\left(8\lambda_2\right)^{-\frac{1}{2p}}$.
\end{lemm}
\noindent
\begin{dem}
The solution that we consider is given by $$f = -\int_0^\infty \left(P_t(g) - \int g\,d\mu\right)dt, $$ where $P_t$ is the semigroup generated by $\cL$. Since $\mu$ satisfies a Poincaré inequality, we can use \eqref{eq:cvLpsemigrp}, and we obtain
\begin{align*}
\int f^p d\mu &= \int \left( \int_0^\infty \left(P_t(g) - \int g\,d\mu\right)dt \right)^p d\mu\\
&\leq \int_0^\infty \left|\left| P_t(g) - \int g\,d\mu \right|\right|_{L^p(\mu)}^p dt\\
&\leq \int_0^\infty K_p^p\, e^{-p\lambda_p\,t} \left|\left| g - \int g\,d\mu \right|\right|_{L^p(\mu)}^p dt\\
&= \frac{K_p^p}{p\lambda_p}\left|\left| g - \int g\,d\mu \right|\right|_{L^p(\mu)}^p .
\end{align*}
Now since $p\geq 2$, we have that $q\in[1,2]$, and therefore $||\nabla f||_{L^q(\mu)}^2 \leq ||\nabla f||_{L^2(\mu)}^2$ and $||\cL f||_{L^q(\mu)} \leq ||\cL f||_{L^p(\mu)}$  so we can write that
\[
||\nabla f||_{L^q(\mu)}^2 \leq -\int f \cL f\,d\mu \leq ||f||_{L^p(\mu)} \left|\left| \cL f \right|\right|_{L^q(\mu)} \leq \frac{K_p}{\sqrt[p]{p\lambda_p}}\left|\left| g - \int g\,d\mu \right|\right|_{L^p(\mu)}^2
\]
which proves \eqref{eq:lemCZforL}.
\end{dem}\newline\newline
Note that in the case $p=2$, Inequality \eqref{eq:lemCZforL} is in fact equivalent to the Poincaré inequality \eqref{def:Poincare}, see e.g. \cite{ane2000inegalites}.

\subsection{Application to $L^p$-Cheeger inequalities under the Poincaré inequality}\label{sec:poincareimpliqueLpCheeger}
In Section \ref{sec:LpCheeger} and as a consequence of the general abstract theorem \ref{thm:main}, we have seen that first order Calderón-Zygmund estimates \eqref{eq:CZgeneral} imply the $L^p$-Cheeger inequalities \eqref{eq:LpCheeger}. Moreover, in Section \ref{sec:CZ}, we have seen that the Poincaré inequality \eqref{def:Poincare} implies the $L^q-L^p$ first order Calderón-Zygmund estimate for all $p\geq 2$, with $q$ the conjugate exponent of $p$. Therefore, we have proven the following theorem, which states that the Poincaré inequality implies the $L^p$-Cheeger inequality for all $p\geq 2$.
\begin{thm}\label{thm:poincareimpliqueCheegerLp}
Let $\mu=e^{-V}dx$, $V\in\C^2(\R^d)$, be a probability distribution satisfying a Poincaré inequality with constant $\lambda_2$. Then for all $p\geq 2$, the following $L^p$-Cheeger inequality holds for all measurable subsets $\Omega\subset\R^d$,
\[
\mu^{+,p}(\partial\Omega) \geq  \frac{1}{\sqrt{2}}\left(8\lambda_2\right)^{\frac{1}{2p}} \, \min\left(\mu(\Omega), 1-\mu(\Omega)\right).
\]
\end{thm}

\subsection{Application to gradient estimates from Cheeger's inequality}\label{sec:reverse}

So far, we have shown how different types of gradient estimates imply different types of isoperimetric inequalities. In that section we state a partial converse: the usual Cheeger inequality implies first order Calderón-Zygmund estimates for all $p\geq 2$.
\begin{thm}\label{thm:reversehierarchy}
Let $\mu=e^{-V}dx$, $V\in\C^2(\R^d)$, be a probability distribution satisfying Cheeger's inequality 
\[
\mu^{+}(\partial\Omega) \geq  h_\mu \, \min\left(\mu(\Omega), 1-\mu(\Omega)\right).
\]
Then for all $p\geq 2$, the following first order $L^q-L^p$ Calderón-Zygmund estimate holds
\begin{equation}\label{eq:CZreversehierarchy}
||\nabla f||_{L^q(\mu)} \leq \sqrt{2}\left(\sqrt{2}h_\mu \right)^{-\frac{1}{p}} ||\cL f||_{L^p(\mu)},
\end{equation}
where $q$ is the conjugate exponent of $p$, and $\cL f = \Delta f -\nabla V\cdot \nabla f$.
\end{thm}
\noindent
\begin{dem}
It is well known, see \cite{cheeger1970lower, yau1975isoperimetric}, that Cheeger's inequality
\[
\mu^{+}(\partial\Omega) \geq  h_\mu \, \min\left(\mu(\Omega), 1-\mu(\Omega)\right),
\]
implies the Poincaré inequality
\[
\lambda_2\,\int \left( f - \int f\,d\mu \right)d\mu \leq \int \left|\nabla f\right|^2 d\mu,
\]
and moreover that the constant $\lambda_2$ in the Poincaré inequality satisfies  $$\lambda_2\geq \frac{1}{4}h_\mu^2. $$ Combined with Lemma \ref{lem:CZforL}, this gives the first order Calderón-Zygmund estimate \eqref{eq:CZreversehierarchy}. 
\end{dem}
\newline \newline
\noindent
\textbf{Aknowledgements.} 
I would like to thank Aldéric Joulin for encouraging me to write this article, and also Xavier Lamy for helping me to tame the notion of total variation. 

\bibliographystyle{plain}
\bibliography{mabibliographie}
\end{document}